\documentclass[12pt]{l4dc2022} 
\title[Stabilization via Discount Policy Gradient]{Learning Stabilizing Controllers of Linear Systems via Discount Policy Gradient}
\usepackage{times}
\usepackage{appendix}
\newtheorem{assum}{Assumption}

\author{\Name{Feiran Zhao} \Email{zhaofr18@mails.tsinghua.edu.cn}\\
  \Name{Xingyun Fu} \Email{fxy20@mails.tsinghua.edu.cn}\\
  \Name{Keyou You} \Email{youky@tsinghua.edu.cn}\\
  \addr Department of Automation and BNRist, Tsinghua University}
\begin{document}	
	\maketitle	

\begin{abstract}
	Stability is one of the most fundamental requirements for systems synthesis. In this paper, we address the stabilization problem for unknown linear systems via policy gradient (PG) methods. We leverage a key feature of PG for Linear Quadratic Regulator (LQR), i.e., it drives the policy away from the boundary of the unstabilizing region along the descent direction, provided with an initial policy with finite cost. To this end, we discount the LQR cost with a factor, by adaptively increasing which gradient leads the policy to the stabilizing set while maintaining a finite cost. Based on the Lyapunov theory, we design an update rule for the discount factor which can be directly computed from data, rendering our method purely model-free. Compared to recent work \citep{perdomo2021stabilizing}, our algorithm allows the policy to be updated only once for each discount factor. Moreover, the number of sampled trajectories and simulation time for gradient descent is significantly reduced to $\mathcal{O}(\log(1/\epsilon))$ for the desired accuracy $\epsilon$. Finally, we conduct simulations on both small-scale and large-scale examples to show the efficiency of our discount PG method.
\end{abstract}


\section{Introduction}

Reinforcement Learning (RL), as a model-free method, has achieved tremendous success empirically in continuous control field~\citep{mnih2015human-level,lillicrap2016continuous}. Instead of identifying an explicit dynamical model, policy optimization methods directly search over the policy space to maximize a performance metric of interest, and have been recognized as an essential approach in modern RL methods. To better understand its performance limits, increasing efforts are devoted to investigating its theoretical guarantees in classical Linear Quadratic Regulator (LQR) problem~\citep{zhou1996robust}, which is a fundamental optimal control framework for stochastic systems. Particularly, the Policy Gradient (PG) methods have been shown in \cite{fazel2018global} to globally converge at a linear rate for the LQR despite its non-convex optimization landscape.

This paper focuses on the stabilization problem for unknown linear stochastic systems, which is among the most fundamental problems in control synthesis, via a PG based approach. The proposed method leverages a key feature of PG for LQR problems, i.e., it drives the policy away from the boundary of unstabilizing region along the descent direction, provided with an initial policy with a finite cost. To this end, we discount the LQR cost with a factor and view it as a free variable. By adaptively updating the discount factor, the gradient descent leads the policy to the stabilizing set while maintaining a finite cost.


\subsection{Related work}
The data-driven stabilization problems for unknown systems have attracted broad attention in both machine learning and control communities. We selectively review some recent literature below.

\textbf{Policy optimization methods.} The optimization landscape of PG for the LQR is firstly made clear by \cite{fazel2018global}, which shows that the PG converges globally due to a gradient dominance property. There are also works focusing on reducing the sample complexity~\citep{malik2019derivative, mohammadi2021convergence}, or solving other important LQR variants e.g., the linear quadratic game~\citep{zhang2019policy}, risk-constrained control \citep{zhao2021primal}, robust control~\citep{gravell2020learn, zhang2021policy}, and linear quadratic Gaussian~\citep{zheng2021sample}. However, all the above works assume that the initial policy is stabilizing, and how to obtain it is posed as an important open problem by \cite{fazel2018global}. Very recently, \cite{perdomo2021stabilizing} takes an initial step towards solving a stabilizing controller via PG methods, which will be discussed in detail later. Though \cite{furieri2020learning, hambly2021policy} consider finite-horizon control and do not require a stabilizing policy, their sample complexity grows linearly in the horizon. 

\textbf{Model-based methods.} The model-based approach solves a stabilizing controller based on an identified model from system trajectories~\citep{abbasi2011regret}. For the case that sampling multiple trajectories is allowed, \cite{dean2020sample} proposes a robust Semi-Definite Program (SDP) to stabilize the uncertain linear system based on the system level synthesis \citep{anderson2019system}. By using a single trajectory, \cite{treven2021learning} and \cite{umenberger2019robust} introduce an ellipsoid region which contains the groundtruth system with high confidence. \cite{chen2021black} considers black-box nonstochastic control and solves a stabilizing controller under adversarial noises via a SDP. 

\textbf{Direct data-driven methods.} There is also a line that directly synthesizes controllers based on a given data set without explicitly identifying a model, originating from the seminal work \citep{willems2005note}. When the data satisfies a Persistent Excitation (PE) condition, \cite{de2019formulas} represents the dynamics using historical trajectories and propose a SDP to stabilize  deterministic linear systems. For the case that data is insufficient for PE condition, \cite{van2020data, van2020noisy} provide both sufficient and necessary for feedback stabilization, and solve a stabilizing controller via a linear matrix inequality for stochastic systems.

\textbf{Discounting methods.} The discounting methods are referred to a class of recently developed system synthesis methods involving damped systems or a varying discount factor~\citep{feng2020escaping, feng2021damping, jing2021learning, lamperski2020computing, perdomo2021stabilizing}. They are firstly investigated in multi-agent control systems to escape local optimal policy~\citep{feng2020escaping,feng2021damping} or compute a stabilizing decentralized controller~\citep{jing2021learning}. Two recent works approach the stabilization problem of centralized linear systems by solving a series of discounted LQR. \cite{lamperski2020computing} considers a deterministic linear system and proposes a model-free policy iteration method with an increasing discount factor to find a stabilizing controller, while the convergence guarantees are totally asymptotic without finite-time rates. A more pertinent work is \cite{perdomo2021stabilizing} which applies policy gradient methods with a damping factor to stabilize both linear and smooth nonlinear systems. It is shown that by using polynomial sampled trajectories, the discount annealing algorithm returns a stabilizing controller in finite iterations. However, all the above work~\citep{feng2020escaping, feng2021damping, jing2021learning, lamperski2020computing, perdomo2021stabilizing} requires a search procedure for the discount factor. For example, \cite{perdomo2021stabilizing} devises a binary search method involving evaluations of noisy function, which complicates the convergence analysis and inevitably increases the sample complexity.

\subsection{Contribution}
We propose a discount policy gradient method with finite-time convergence guarantees to find a stabilizing controller for linear systems. Particularly, our method alternatively updates the policy and discount factor by solely using system trajectories from a simulator. In sharp contrast to \cite{perdomo2021stabilizing}, our method has at least three major merits. First, we use Lyapunov theory to design an explicit update rule for the discount factor that can be directly computed from data instead of invoking a search procedure, which also yields much simpler convergence analysis. Second, thanks to the Lyapunov argument the policy can be updated by one-step gradient descent each iteration, while \cite{perdomo2021stabilizing} requires the policy to converge in principle. Third, we apply new analysis techniques for the policy gradient step \citep{mohammadi2021convergence}, which only requires the simulation time and total number of sampled trajectories to be linear in the desired accuracy. As a comparison, polynomial sample complexity is required in \cite{perdomo2021stabilizing}.

\textbf{Notations.} We use $\rho(\cdot)$ to denote the spectral radius of a matrix, and $\|\cdot\|$ to denote the $2$-norm. Let $\underline{\sigma}(\cdot)$ be the minimal eigenvalue of a matrix. $\text{Tr}(\cdot)$ denotes the trace function. Let $S^{d-1} \subset \mathbb{R}^d$ be the unit sphere of dimension $d-1$. We use $\mathcal{O}(\epsilon)$ to denote some constant proportional to $\epsilon$. 

\section{Problem Formulation}

Consider the following discrete-time linear time-invariant system\footnotemark[1]
\begin{equation}\label{equ:sys}
x_{t+1} = Ax_t + Bu_t, ~~x_0 \sim \mathcal{D},
\end{equation}
where $x_t \in \mathbb{R}^n$ is the state, $u_t\in \mathbb{R}^m$ is the control input. The matrices $A\in \mathbb{R}^{n \times n}$ and $B \in \mathbb{R}^{n \times m}$ are the unknown model parameters. The initial state $x_0$ is sampled from a distribution $\mathcal{D}$, on which we make the following mild assumption\footnotemark[2].

\footnotetext[1]{Though we consider the randomness stemming from the initial state distribution, the noisy case with $x_{t+1} = Ax_t + Bu_t + w_t$ can be addressed using the same method. Please refer to their connections discussed in \cite{malik2019derivative}.}
\begin{assum}
	\label{assumption:D}
	 The distribution $\mathcal{D}$ has zero mean and unit covariance $\mathbb{E}[w_tw_t^{\top}] = I$. Moreover, its support is bounded by $\|x_0\|\leq d$ with a constant $d>0$. 
\end{assum}
\footnotetext[2]{The results in this paper also hold for distribution $\mathcal{D}$ with bounded sub-Gaussian norm, as verified in our experiments; see \cite{mohammadi2021convergence}.}

In this paper, we focus on the stabilization problem of (\ref{equ:sys}) via linear state feedback $u(x_t) = -Kx_t$. Clearly, a minimal assumption is that there exists a gain $K$ such that $\rho(A-BK)<1$.
\begin{assum}
	\label{assumption}
	The pair $(A,B)$ is stabilizable.
\end{assum}

When the explicit model $(A,B)$ is unknown, the PG is widely studied for solving LQR problems, which directly searches over the feedback gain matrix space to solve
\begin{equation}\label{prob:dlqr}
\text { minimize } ~ J_{\gamma}(K) := \mathbb{E}_{x_0} \sum_{t=0}^{\infty}\gamma^{t} (x_{t}^{\top} Q x_{t}+u_{t}^{\top} R u_{t})~~~~
\text {subject to} ~(\ref{equ:sys}),  \text{and}~u_t = -K x_t
\end{equation}
using gradient methods, where $0<\gamma \leq 1$ is a discount factor. However, it requires an initial policy $K^0$ to render a finite cost, i.e., $\sqrt{\gamma}\rho(A-BK^0)<1$, which is non-trivial in the absence of an explicit model. In fact, it is an open problem posed by \cite{fazel2018global} to find an initial stabilizing controller for the case $\gamma = 1$, which is however a common basic assumption in other works on PG for LQR problems \citep{zheng2021sample,zhang2021policy,zhao2021primal}.

\begin{figure}[t]
	\centering
	\includegraphics[height=30mm]{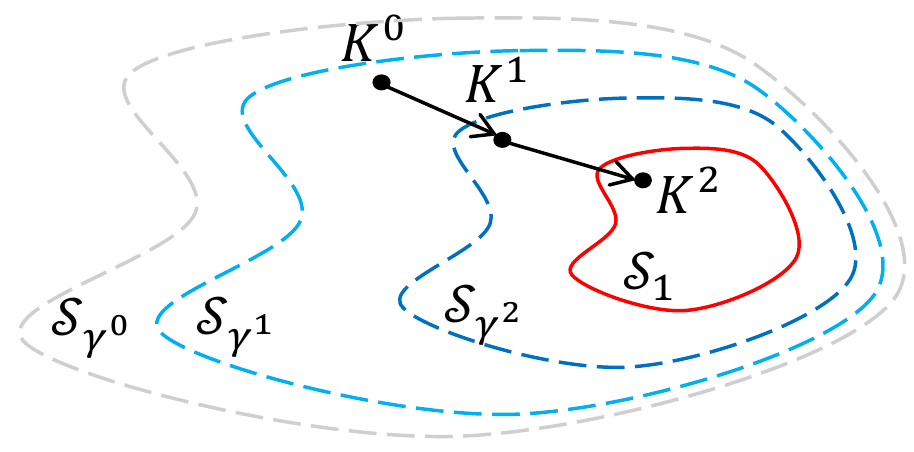}
	\caption{An illustration of the optimization landscape. Let $\mathcal{S}_{\gamma} = \{K | \sqrt{\gamma}\rho(A-BK) < 1 \}$. Our method starts from a policy $K^0 =0 \in \mathcal{S}_{\gamma^0}$ with $\gamma^0 < \rho(A^2)$. We first compute a discount factor $\gamma^1 > \gamma^0$ such that $K^0 \in \mathcal{S}_{\gamma^1}$. Then, $J_{\gamma^1}(K)$ is optimized via gradient methods to render the updated policy $K^1$ away from the boundary of $\mathcal{S}_{\gamma^1}$. This process is iterated until $\gamma^i \geq 1$, which implies that $K^{i-1}$ is stabilizing for (\ref{equ:sys}).}
	\label{pic:pgland}
\end{figure}

In this paper, we view $\gamma$ as a variable and propose a discount policy gradient method to stabilize (\ref{equ:sys}), which alternatively updates $K$ and $\gamma$ by solely using data from a simulator. In particular, we use gradient descent to update $K$ towards the stabilizing region, and ensure $\sqrt{\gamma}\rho(A-BK)<1$ by designing an update rule for $\gamma$; see Fig. \ref{pic:pgland} for an illustration. Moreover, the rule can be directly computed by the data, and hence a search procedure is not required.

%

\section{Policy gradient and Lyapunov theory}

In this section, we first provide preliminaries on the discounted LQR problem. Then, we leverage Lyapunov theory to design an update rule for the discount factor.

Consider the discounted LQR problem in (\ref{prob:dlqr}). Let $Q>0$ and $R>0$ be user-specified matrices. The following lemma provides a condition for $J_{\gamma}(K)$  to be finite.
\begin{lemma}\label{lem:closed}
	$J_{\gamma}(K) < \infty$ if and only if $\sqrt{\gamma}\rho(A-BK)<1$. Moreover, if $\sqrt{\gamma}\rho(A-BK)<1$, then the cost has a closed-form expression 
	$J_{\gamma}(K) = \text{Tr}(P_K^{\gamma}),$
	where $P_K^{\gamma}$ is a unique positive definite solution to the Lyapunov equation 
	$
	P_K^{\gamma} = Q + K^{\top}RK + \gamma (A-BK)^{\top}P_K^{\gamma}(A-BK).
	$
\end{lemma}

Given a policy $K$ with $J_{\gamma}(K) < \infty$, the policy gradient method updates $K$ by 
\begin{equation}\label{equ:pg}
	K^{+} = K - \eta \hat{\nabla} J_{\gamma}(K),
\end{equation}
where $\hat{\nabla} J_{\gamma}(K)$ is an estimation of the gradient with respect to $K$, and $\eta$ denotes an appropriate step size. Under proper gradient estimate methods such as minibatching~\citep{fazel2018global} and one- or two-point estimation \citep{malik2019derivative,mohammadi2021convergence}, it has been shown to converge linearly to an optimal policy while maintaining a finite cost. That is, the gradient descent (\ref{equ:pg}) drives the policy away from the boundary of $\mathcal{S}_{\gamma}$. Hence, we can find a larger $\gamma^{+}>\gamma$ such that $K^{+} \in \mathcal{S}_{\gamma^{+}}$. Next, we introduce the Lyapunov stability theory to design such a discount factor.

The following well-known result \citep{vidyasagar2002nonlinear} provides a sufficient condition for the stability of linear systems. 
\begin{lemma}\label{lem:Lyapunov}
	Consider system (\ref{equ:sys}). Suppose that there exists a function $V$: $\mathbb{R}^n \rightarrow \mathbb{R}$ continuous at the origin, finite for $x\in \mathbb{R}^n$, and such that 
	\begin{equation}\label{def:Lyapnov}
		V(0)=0 , V(x)>0, \forall x \neq 0, ~~\text{and}~~ V(x_{t+1}) - V(x_t) < 0.
	\end{equation}
	Then, (\ref{equ:sys}) is asymptotically stable at $x=0$.
\end{lemma}

\begin{definition}
	A function $V(x)$ satisfying (\ref{def:Lyapnov}) is called a Lyapunov function. 
\end{definition}

Before proceeding, we note that $\sqrt{\gamma}\rho(A-BK)<1$ if and only if the scaled dynamical system
\begin{equation}\label{equ:rescale_dyna}
	x_{t+1} = \sqrt{\gamma}(A - BK) x_t
\end{equation}
is stabilizing, which enables us to find a feasible $\gamma$ by designing a Lyapunov function for (\ref{equ:rescale_dyna}).

\begin{lemma}\label{theo:rule}
	Suppose that $\sqrt{\gamma}\rho(A-BK)<1$. Let $P>0$ be the solution of the Lyapunov function
	\begin{equation}\label{equ:Pi}
		P = Q + K^{\top}RK + \gamma (A-BK)^{\top}P(A-BK).
	\end{equation}	
	Then, $J_{\gamma'}(K)<+\infty$ if $\gamma'$ satisfies
	\begin{equation}\label{equ:rule}
	(1-\gamma/\gamma')P < Q + K^{\top}RK.
	\end{equation}
\end{lemma}

\begin{proof}
	To prove $J_{\gamma'}(K)<+\infty$, it is equivalent to show $\sqrt{\gamma'}\rho((A-BK))<1$. Consider the following autonomous system
	\begin{equation}\label{sys:lyapnov}
	x_{t+1} = \sqrt{\gamma'}(A-BK)x_t.
	\end{equation}	
	We show that under the given condition (\ref{equ:rule}), $V(x) = x^{\top}Px$ with $P$ defined in (\ref{equ:Pi}) is a Lyapunov function of the scaled system (\ref{sys:lyapnov}). Then, by Lemma \ref{lem:Lyapunov}, (\ref{sys:lyapnov}) is asymptotically stable at $x = 0$.
	
	Clearly, $V(x)$ satisfies that $V(0) = 0$ and $V(x)>0,\forall x\neq 0$. Further note that  
	\begin{align*}
	V(x_{t+1}) - V(x_t) &= \gamma'x_t^{\top}(A-BK)^{\top}P(A-BK)^{\top}x_t - x_t^{\top}P x_t \\
	&= x_t^{\top}(\frac{\gamma'}{\gamma}(P-Q-K^{\top}RK) - P) x_t,
	\end{align*}
	where the last equation follows from (\ref{equ:Pi}). By using the condition in (\ref{equ:rule}), it follows that $V(x_{t+1}) - V(x_t)<0$, which proves that $V(x) = x^{\top}Px$ is indeed a Lyapunov function for (\ref{sys:lyapnov}).
\end{proof}

Lemma \ref{theo:rule} provides a sufficient condition on $\gamma'$ to ensure $J_{\gamma'}(K)<+\infty$, provided with $\sqrt{\gamma}\rho(A-BK)<1$. However, it is impossible to evaluate $\gamma'$ as the computation of $P$ involves the unknown model parameters $(A, B)$. In the following theorem, we show that by utilizing the closed-form expression of the cost in Lemma \ref{lem:closed}, the condition (\ref{equ:rule}) can be approximately computed using data from a simulator. 

\begin{theorem}\label{coro}
Suppose that $\sqrt{\gamma}\rho(A-BK)<1$. Then, $J_{\gamma'}(K)<+\infty$ if 
	\begin{equation}\label{equ:coro}
		\gamma' \leq (1 + \frac{\underline{\sigma}(Q+K^{\top}RK)}{J_{\gamma}(K) - \underline{\sigma}(Q+K^{\top}RK)}) \gamma.
	\end{equation}
\end{theorem}

\begin{proof}
	To ensure (\ref{equ:rule}), it suffices for $\gamma'$ to satisfy
	$
	1-{\gamma}/{\gamma'} < {\underline{\sigma}(Q+K^{\top}RK)}/{\|P\|}.
	$
	Noting that
	$
	J_{\gamma}(K) =\text{Tr}(P) \geq \|P\|,
	$	the proof is completed.
\end{proof}

Since the parameters $(Q, R)$ are user-defined, the term $\underline{\sigma}(Q+K^{\top}RK)$ can be directly computed. Moreover, the cost $J_{\gamma}(K)$ can be evaluated by sampling system trajectories from the simulator. In the sequel, we propose our algorithm based on Theorem \ref{coro} and gradient descent (\ref{equ:pg}).

\begin{algorithm2e}[t]
	\caption{The discount policy gradient algorithm}
	\label{alg:InfPG}	
	\LinesNumbered
	\KwIn{Initial policy $K^0 = 0$ and discount factor $\gamma^0$, simulation time $\tau$, number of trajectories $N$ for cost evaluation.}
	\For {$i=0,1,\cdots$}
	{
		Evaluate $\hat{J}_{\gamma^i}^{\tau}(K^{i}) = \frac{1}{N}\sum_{j=0}^{N-1} V_{\gamma^i}^{\tau}(K^i, x_0^j) $ with $x_0^j$ sampled independently from $\mathcal{D}$\;
		Compute the discount factor $\gamma^{i+1} = (1+\alpha^i)\gamma^i$ with $\alpha^i$ given by \
		\begin{equation}\label{equ:dis_update}
		\alpha^i = \underline{\sigma}(Q+(K^i)^{\top}RK^i)/(2\hat{J}_{\gamma^i}^{\tau}(K^{i}) - \underline{\sigma}(Q+(K^i)^{\top}RK^i));
		\end{equation} 
		\If{$\gamma^{i+1} \geq  1 $}{Return a stabilizing policy $K^i$ \;}
		Update $K^{i+1}$ via policy gradient in (\ref{equ:pg}) starting from policy $K^i$, such that 
		\begin{equation}\label{equ:Jbar}
			J_{\gamma^{i+1}}(K^{i+1}) < \bar{J}.
		\end{equation}
	}
\end{algorithm2e}

\section{Main Results}
In this section, we present the discount policy gradient algorithm to find a stabilizing controller for (\ref{equ:sys}), the convergence of which is shown to be linear with respect to the initial $\gamma^0$.

The algorithm is detailed in Algorithm \ref{alg:InfPG}. The initial discount factor should be selected to satisfy $\sqrt{\gamma^0}\rho(A)<1$ such that the Lyapunov equation (\ref{equ:Pi}) has a solution. Line 2 of Algorithm \ref{alg:InfPG} evaluates the cost by Monte Carlo sampling with truncated cost 
$$V_{\gamma}^{\tau}(K,x_0) = \sum_{t=0}^{\tau-1}\gamma^{t} (x_{t}^{\top} Q x_{t}+u_{t}^{\top} R u_{t}),~~\text {subject to} ~(\ref{equ:sys}), x_0 \in \mathcal{D}, \text{and}~u_t = -K x_t., $$
 where $\tau$ is the simulation time, and a realization of the initial state is sampled from distribution $\mathcal{D}$. In line 3, we apply Theorem \ref{coro} to compute an update rate for the discount factor. Note that we have considered the estimation error of $J_{\gamma^i}(K^i)$ in line 2 in the denominator of (\ref{equ:dis_update}). In line 6, we only require $\bar{J}$ to be a positive constant satisfying $\bar{J} > J^*_{1}$ (to be specified later) for the convergence analysis, which can be achieved by one-step gradient descent in practice.

First, we show that under the condition $J_{\gamma}(K) \leq \bar{J}$ in (\ref{equ:Jbar}), the estimation error $|\hat{J}_{\gamma}^{\tau}(K)-J_{\gamma}(K)|$ induced by Monte Carlo sampling and finite simulation time can be well controlled with a large probability. Different from \citet[Lemma 26]{fazel2018global}, we propose a new proof technique to show that the simulation time for function evaluation to achieve $\epsilon$-accuracy is proportional to only $\log(1/\epsilon)$, hence improved from $\text{poly}(1/\epsilon)$ in \cite{perdomo2021stabilizing}. For the sake of exposition, we omit problem-dependent constants in $\mathcal{O}(\cdot)$. 

\begin{lemma}\label{lem:error}
	For a given constant $0<\delta<1$, let the simulation time be $\tau = \mathcal{O}(\log(1/J_{\gamma}(K)))$ and the number of samples be $N = \mathcal{O}(\frac{1}{J^2_{\gamma}(K)}\log(1/\delta))$. Then, with at least probability $1-\delta$, it holds
	\begin{equation}\label{equ:eval_bound}
		|\hat{J}_{\gamma}^{\tau}(K) - J_{\gamma}(K)| \leq \frac{1}{2}J_{\gamma}(K).
	\end{equation}
	
\end{lemma}

\begin{proof}
	We first establish an upper bound of the bias $V_{\gamma}(K, x_0) - V_{\gamma}^{\tau}(K, x_0)$ induced by finite simulation time $\tau$, which has exponential dependence on $\tau$. For a policy with $\sqrt{\gamma}\rho(A-BK)<1$, 
	\begin{align*}
		&V_{\gamma}(K, x_0) - V_{\gamma}^{\tau}(K, x_0) \\
		& =  \sum_{t=0}^{\infty}\gamma^{t} (x_{t}^{\top} Q x_{t}+u_{t}^{\top} R u_{t}) -\sum_{t=0}^{\tau -1}\gamma^{t} (x_{t}^{\top} Q x_{t}+u_{t}^{\top} R u_{t})  
		 = \sum_{t=\tau}^{\infty}\gamma^{t} (x_{t}^{\top} Q x_{t}+u_{t}^{\top} R u_{t}) \\
		 &\leq \text{Tr}(P)\|\sqrt{\gamma}(A-BK)\|^{2\tau}\|x_0\|^2 
		 \leq \bar{J}d^2\cdot\|\sqrt{\gamma}(A-BK)\|^{2\tau},
	\end{align*}
	where the last inequality follows from $\text{Tr}(P) = J_{\gamma}(K) < \bar{J}$ and $\|x_0\|\leq d$. The norm square $\|\sqrt{\gamma}(A-BK)\|^2$ can be bounded by
	$$
	\|\sqrt{\gamma}(A-BK)\|^2 \leq \sup \frac{\|\sum_{t=0}^{\infty}\gamma^t((A-BK)^t)^{\top}X(A-BK)^t\|}{\|X\|}\leq J_{\gamma}(K)/\underline{\sigma}(Q) \leq \bar{J}/\underline{\sigma}(Q),
	$$
	where the second inequality follows from \cite[Lemma 17]{fazel2018global}. Thus, it follows that
	$$
	V_{\gamma}(K, x_0) - V_{\gamma}^{\tau}(K, x_0) \leq \bar{J}d^2(\bar{J}/\underline{\sigma}(Q))^{\tau}.
	$$
	
	Then, we use concentration inequalities to bound the total error $|J_{\gamma}(K) - \hat{J}_{\gamma}^{\tau}(K)|$. Let $x_0^j, j\in \{1,2,\dots,N \}$ be $N$ random initial state sampled independently from $\mathcal{D}$. 
	Since the support of $\mathcal{D}$ is bounded by $\|x_0\| \leq d$, the random variable $V_{\gamma}(K, x_0^j)$ is bounded by $ 0\leq V_{\gamma}(K, x_0^j) = \text{Tr}\{P x_0x_0^{\top}\} \leq \bar{J}d^2, \forall j \in \{1,2,\dots,N \} $. For a give constant $\epsilon$, we let the simulation time be $\tau = \frac{-\log(2\bar{J}d^2/\epsilon)}{\log(\bar{J}/\underline{\sigma}(Q))}$ such that $V_{\gamma}(K, x_0^j) - V_{\gamma}^{\tau}(K, x_0^j) \leq \frac{1}{2}\epsilon$, which implies that 
	\begin{equation}\label{equ:error}
		|\frac{1}{N}\sum_{j=0}^{N-1} V_{\gamma}^{\tau}(K, x_0) -\frac{1}{N}\sum_{j=0}^{N-1} V_{\gamma}(K, x_0)| \leq \frac{1}{2}\epsilon.
	\end{equation}
	Hence, the Hoeffding's inequality yields that
	\begin{align*}
		\text{Pr}(|\hat{J}_{\gamma}^{\tau}(K) - J_{\gamma}(K)| \leq \epsilon) &= \text{Pr}(|\frac{1}{N}\sum_{j=0}^{N-1} V_{\gamma}^{\tau}(K, x_0) - \mathbb{E}_{x_0}V_{\gamma}(K,x_0)| \leq \epsilon) \\
		&\geq \text{Pr}(|\frac{1}{N}\sum_{j=0}^{N-1} V_{\gamma}(K, x_0) - \mathbb{E}_{x_0}V_{\gamma}(K,x_0)| \leq \frac{\epsilon}{2}) \\
		& \geq 1- 2\exp(-\frac{N\epsilon^2}{2\bar{J}d^2}),
	\end{align*}
	where the first inequality follows from (\ref{equ:error}). Let $\delta = 2\exp(-\frac{2N\epsilon^2}{\bar{J}d^2})$. Then, we conclude that at least with probability $1-\delta$, the estimation error is bounded by $|\hat{J}_{\gamma}^{\tau}(K) - J_{\gamma}(K)| \leq \epsilon$. Letting $\epsilon = J_{\gamma}(K)/2$ completes the proof.
\end{proof}

By lemma \ref{lem:error}, we have that with large probability, $J_{\gamma}(K) < 2\hat{J}_{\gamma}^{\tau}(K)$. Thus, it suffices to choose $\alpha^i$ in (\ref{equ:dis_update}) such that Theorem \ref{coro} holds. To yield fast convergence of the discount factor, the update rate $\alpha^i$ must be lower bounded by a positive constant. Hence, we require that after the gradient descent in line 6, the cost is upper bounded by a positive constant $J_{\gamma}(K) < \bar{J}$. The following lemma provides a uniform lower bound for $\bar{J}$.
\begin{lemma}\label{lem:barJ}
	For $0<\gamma_1<\gamma_2\leq 1$, it follows that 
	$
		J^*_{\gamma_1} < J^*_{\gamma_2},
	$
	where $J^*_{\gamma}$ denotes the optimal value of $J_{\gamma}(K)$, i.e., $J^*_{\gamma} < J_{\gamma}(K), \forall K \in \mathcal{S}_{\gamma}$.
\end{lemma}

Thus, $\bar{J}$ can be any positive constant larger than $J^*_1$. In fact, a larger $\bar{J}$ implies less gradient descent steps each iteration and more total iterations of Algorithm \ref{alg:InfPG}, as to be shown later.

Then, we establish the convergence of the gradient descent in line 6. We apply a two-point gradient estimation for $\hat{\nabla}J(K)$, as it yields better sample complexity than the one-point setting. Motivated by \cite{mohammadi2020linear}, we show that with a large probability, (\ref{equ:pg}) converges linearly using only $\log(1/\epsilon)$ simulation time and total sampled trajectories with a desired accuracy $\epsilon$ of the cost, which significantly improves the polynomial sample complexity in \cite{perdomo2021stabilizing}. By letting  $\epsilon = \bar{J} - J_{\gamma^{i+1}}^*$, we have the following result.
\begin{algorithm2e}[t]
	\caption{Two-point gradient estimation}
	\label{alg:gradient}	
	\KwIn{Policy $K$, distribution $\mathcal{D}$, discount factor $\gamma$, smoothing radius $r$, simulation time $\tau$, number of random samples $M$.}
	\For{$j=1,2,\cdots,M$}
	{
		Sample a perturbation matrix $U_{j}$ uniformly from the sphere $\sqrt{mn}S^{mn-1}$ \;
		Set $K_{j,1} = K + rU_{j}$ and $K_{j,2} = K - rU_{j}$ \;
		Sample an initial state $x_0^j$ from distribution $\mathcal{D}$ \;
		Simulate system (\ref{equ:sys}) to compute $V_{\gamma}^{\tau}(K_{j,1}, x_0^j)$ and $V_{\gamma}^{\tau}(K_{j,2}, x_0^j)$ \;
	}
	\KwOut{Gradient estimation $\hat{\nabla}J_{\gamma}(K) = \frac{1}{2rM} \sum_{j=1}^{M}(V_{\gamma}^{\tau}(K_{j,1}, x_0^j) - V_{\gamma}^{\tau}(K_{j,2}, x_0^j))U_j.$}
\end{algorithm2e}
\begin{lemma}\label{lem:gradient}
	Consider the policy gradient method (\ref{equ:pg}) at $i$-th iteration in Algorithm \ref{alg:InfPG} with two-point gradient estimation by Algorithm \ref{alg:gradient}. Let the simulation time $\tau$ and number of trajectories $M$ satisfy
	$
	\tau = \mathcal{O}(\log({1}/{ \bar{J} - J_{\gamma^{i+1}}^*}))$ and $ M = \mathcal{O}(1).
	$ Then, for a smoothing radius $r<\mathcal{O}(\sqrt{ \bar{J} - J_{\gamma^{i+1}}^*})$ and some constant step size $\eta$, conducting (\ref{equ:pg}) from initial policy $K^i$ in $\mathcal{O}(\log(J_{\gamma^{i+1}}(K^{i})/( \bar{J} - J_{\gamma^{i+1}}^*)))$ iterations achieves $J_{\gamma^{i+1}}(K^{i+1}) < \bar{J}$ with large probability.
\end{lemma}

\begin{theorem}\label{the:finite-time}
	In at most $\frac{3\bar{J}-\underline{\sigma}(Q)}{\underline{\sigma}(Q)} \log \frac{1}{\gamma^0}$ iterations, Algorithm \ref{alg:InfPG} returns a stabilizing controller of (\ref{equ:sys}). 
\end{theorem}

\begin{proof}
	At $i$-th iteration in Algorithm \ref{alg:InfPG}, suppose that the condition $J_{\gamma^{i}}(K^{i}) < \bar{J}$ is satisfied. Then, the update rate $\alpha^i$ is uniformly lower bounded by 
	$$
	\alpha^i = \frac{\underline{\sigma}(Q+(K^i)^{\top}RK^i)}{2\hat{J}_{\gamma^i}^{\tau}(K^{i}) - \underline{\sigma}(Q+(K^i)^{\top}RK^i)} \geq  \frac{\underline{\sigma}(Q)}{3J_{\gamma^i}(K^i) - \underline{\sigma}(Q)} \geq \frac{\underline{\sigma}(Q)}{3\bar{J} - \underline{\sigma}(Q)},
	$$
	where the first inequality follows from (\ref{equ:eval_bound}).
	
	Hence, we conclude that the number of iterations for Algorithm \ref{alg:InfPG} is no more than $\log(1/\gamma^0)/\log(1+{\underline{\sigma}(Q)}/{(3\bar{J} - \underline{\sigma}(Q))})$. Then, simplifying it via $\log(1+x)\approx x$ completes the proof.
\end{proof}

As a comparison, \cite{perdomo2021stabilizing} requires the number of iterations for $\gamma$ to be $64 (J_{1}^*)^4 \log \frac{1}{\gamma^0}$ by letting $\underline{\sigma}(Q) = 1$, which is a $4$-th order polynomial of ours.

\section{Experimental analysis}
\begin{figure}[t]
	\centering
	\includegraphics[height=60mm]{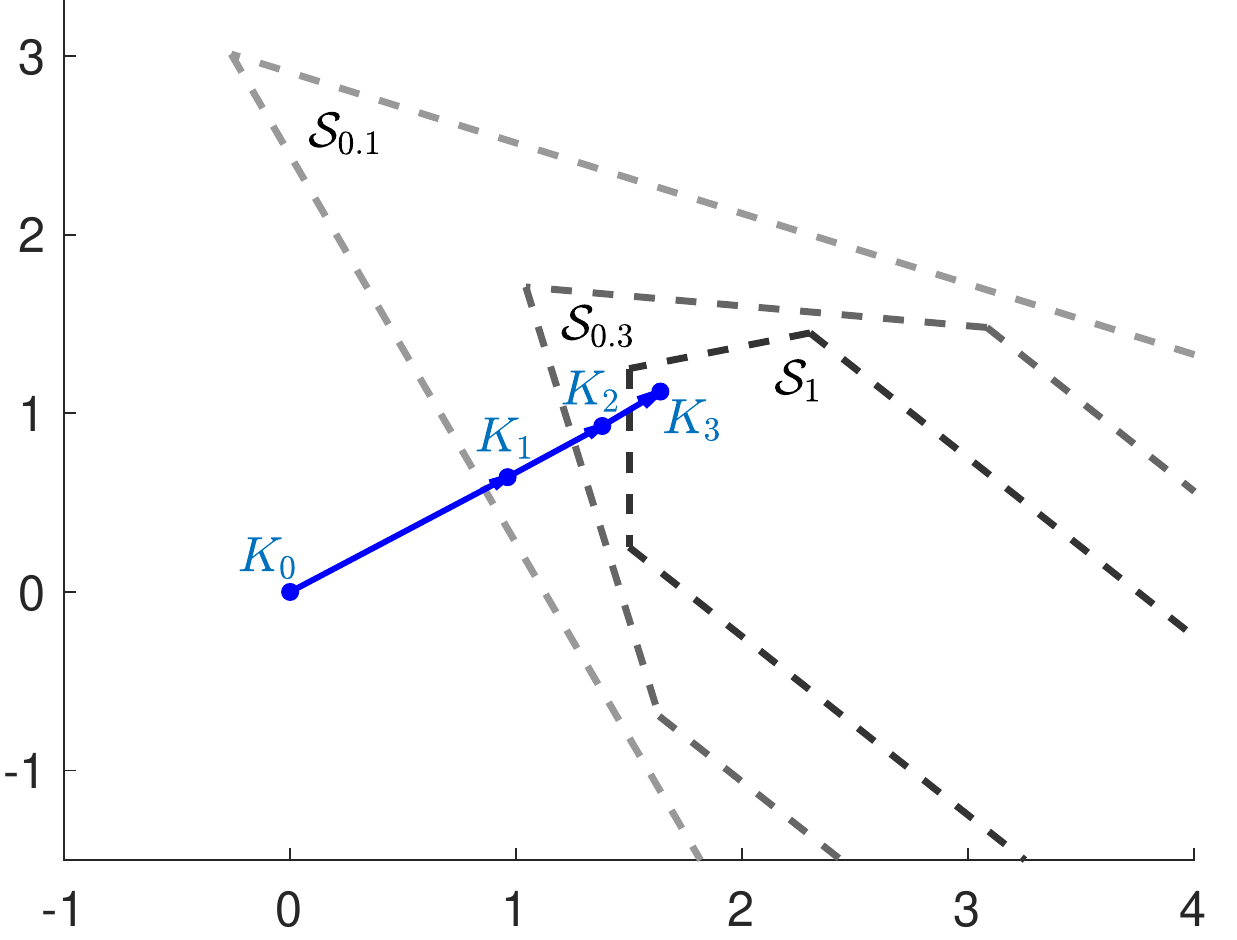}
	\caption{The optimization process of Algorithm \ref{alg:InfPG} on the small-scale example. We plot the boundaries of the sets $\mathcal{S}_{\gamma} = \{K | \sqrt{\gamma}\rho(A-BK) < 1 \}$ with $\gamma = 0.1, 0.3, 1$, and select three intermediate points denoted by $K_1,K_2,K_3$ for illustration. Each iteration of Algorithm \ref{alg:InfPG} drives the policy one step towards the stabilizing set.
	}
	\label{pic:policy}
\end{figure}

This section verifies the effectiveness of the proposed discount policy gradient methods in Algorithm \ref{alg:InfPG}. The simulation is carried out using MATLAB 2021b on a laptop with 2.8GHz CPU. The code is provided in \url{https://github.com/fuxy16/Stabilize-via-PG}.

We first conduct experiments on a two-dimensional example for illustration. Consider the following unstable dynamical model with single control input and penalty matrices
$$
A = \begin{bmatrix}
4 &3 \\
3 &1.5
\end{bmatrix}, ~~ B = \begin{bmatrix}
2 \\
2
\end{bmatrix}, ~~Q = \begin{bmatrix}
1 \\
1
\end{bmatrix}, ~~R = 2.
$$
Clearly, $(A,B)$ is controllable. Let the initial state distribution $\mathcal{D}$ be the standard normal distribution. The number of sampled trajectories in function evaluation (line 2 in Algorithm \ref{alg:InfPG}) is set to $N = 50$, and the simulation time $\tau$ is set to $\tau = 100$. We select an initial policy $K^0 = 0$ and discount factor $\gamma^0 = 10^{-3} < 1/\rho^2(A) = 1/36$. We apply one-step gradient descent each iteration with a constant step size $\eta = 10^{-3}$, where the smooth radius and the number of trajectories for  gradient estimation in Algorithm \ref{alg:gradient} is $r = 2\times 10^{-3}$ and $M =10$. Fig. \ref{pic:policy} illustrates the optimization process of Algorithm \ref{alg:InfPG}. In less than $250$ iterations, Algorithm \ref{alg:InfPG} returns a stabilizing controller. 
\begin{figure}[t]
	\centering
	\subfigure[Iteration of Algorithm \ref{alg:InfPG}.]{
		\includegraphics[height=45mm]{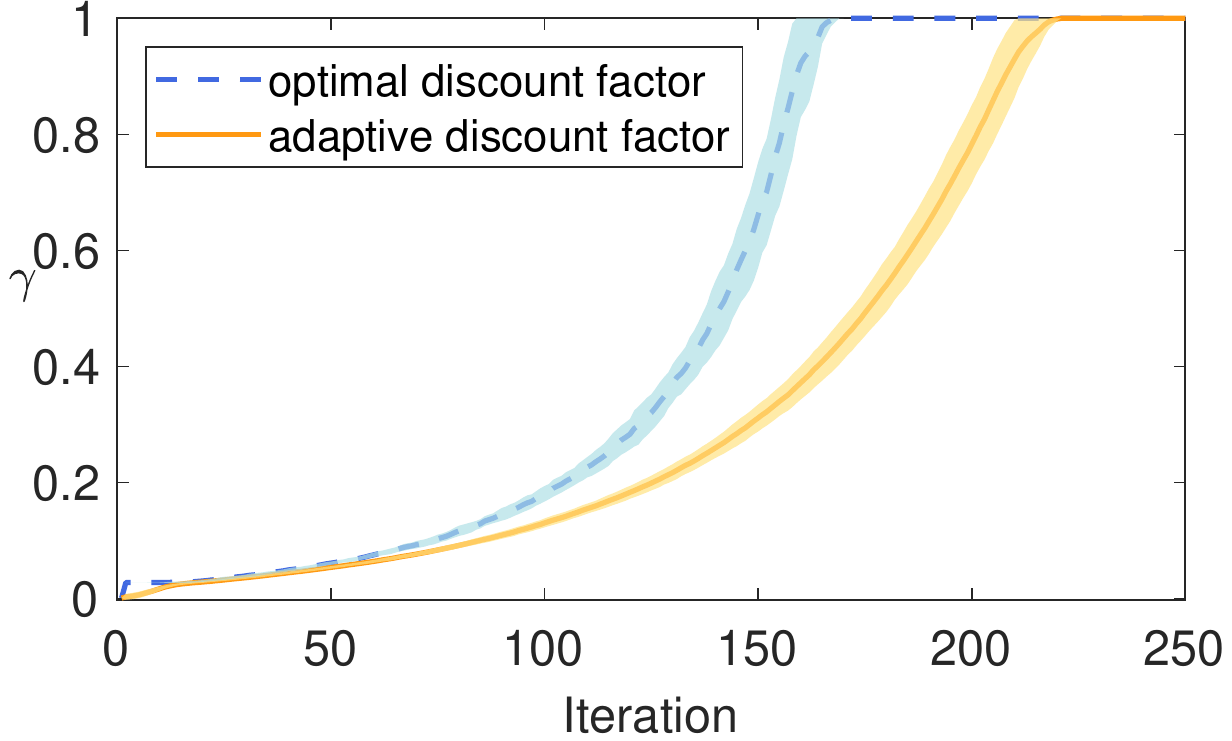}}
	\subfigure[Model-based implementation.]{
		\includegraphics[height=45mm]{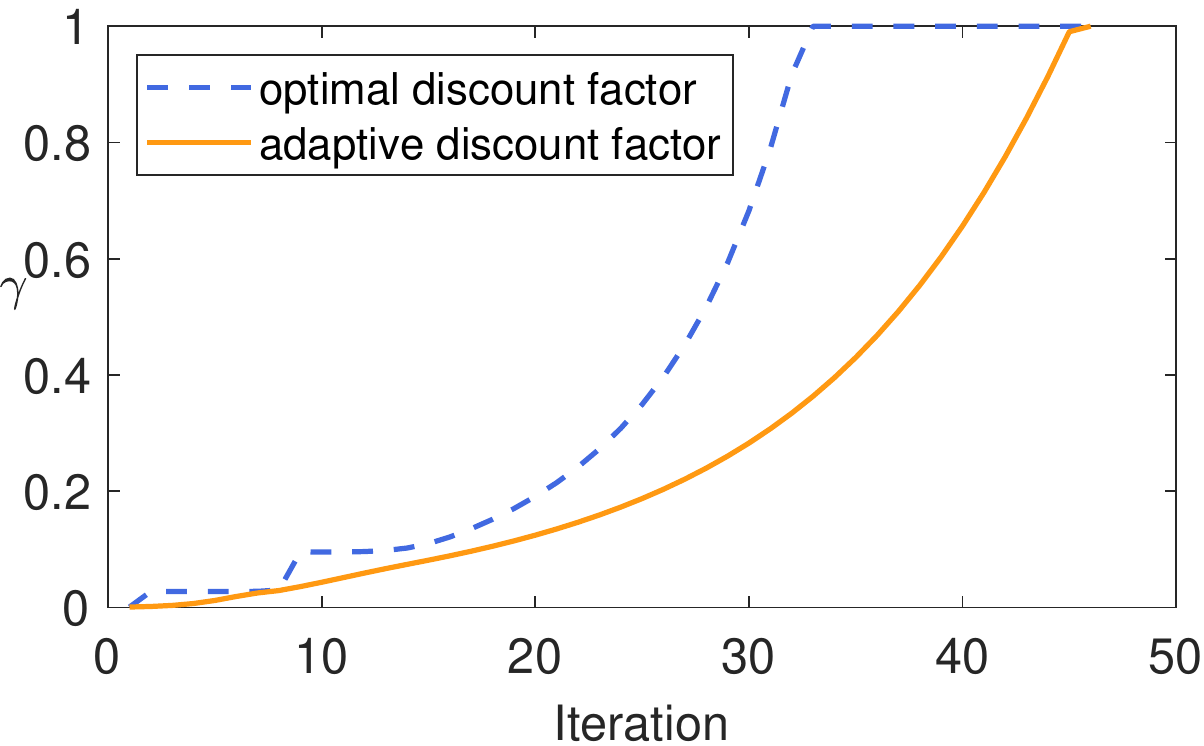}}
	\caption{Convergence of the discount factor. The centreline denotes the mean of 20 independent trials and the shaded region demonstrates their standard deviation.}
	\label{pic:discount}
\end{figure}

Next, we show the convergence of the discount factor in Fig \ref{pic:discount}, where we report the results of $20$ independent trails. From $\gamma^0 = 10^{-3}$, the adaptive discount factor (yellow solid line) grows almost exponentially to $1$ within $250$ iterations. The sample complexity can be calculated by $250\times(M+N)$, i.e.,  total number of $1.5\times 10^4$ sampled trajectories with simulation time $\tau = 100$. For comparison, we also plot the maximal discount factor under $K^i$ each iteration, given as $\gamma^i_{\text{opt}} = 1/\rho^2(A-BK^i)$, denoted by the ``optimal'' discount factor (blue dashed line). We observe that the update rate $\alpha^i$ in (\ref{equ:dis_update}) approximates the upper bound well in the first $100$ iterations, and the gap increases reasonably due to the approximation of condition (\ref{equ:rule}). Moreover, the variance induced by independent trials is competitively small, considering our low sample complexity. We also display a model-based implementation of Algorithm \ref{alg:InfPG}, where we assume that $(A,B)$ are known and compute the update rate $\alpha$ by (\ref{equ:coro}) in Theorem \ref{coro}. It is shown that with an accurate evaluation $J_{\gamma}(K)$ and the policy gradient, the required iterations can be reduced to less than $50$. 

Finally, we show the efficiency of Algorithm \ref{alg:InfPG} on large-scale problems. We randomly sample $20$ independent pairs of system matrices $A,B \in \mathbb{R}^{100\times 100}$ where each element of $B$ is subject to a standard normal distribution, and each element of $A$ is sampled from a normal distribution with variance $0.01$ to ensure the controllability. We conduct $20$ independent trials using the same parameter setting as the two-dimensional example. The mean of the number of iterations is $1750$, which implies a total number $2\times 10^{5}$ of sampled trajectories, showing the efficiency of our method.

\section{Conclusion and future work}

This paper proposes a flexible policy gradient based framework to find a stabilizing controller for an unknown linear system by adaptively updating the discount factor. In particular, our algorithm only requires $\mathcal{O}(\log(1/\gamma^0))$ iterations, and both the simulation time as well as the number of sampled trajectories for gradient descent grow linearly in the desired accuracy, which significantly improves from the existing methods that requires polynomial sample complexity~\citep{perdomo2021stabilizing}. 

We now discuss some possible future directions. The first one is to replace the policy gradient step in Algorithm \ref{alg:InfPG} with other learning-based methods, such as Q-learning. This is because we only requires the cost to decrease at each iteration, which can also be done by other reinforcement learning methods. A second one is to study policy gradient of finite-horizon LQR for system stabilization with a single feedback gain $K$. The idea is to optimize $K$ with an increasing horizon. The main challenge is that an explicit expression for $K$ might not exist and the optimization landscape is unclear yet. The third direction is to further investigate the role of the discount factor in the policy optimization of LQR problems. Existing literature in the RL field has shown that the discount factor works as a regularizer for generalization~\citep{amit2020discount}, and can even accelerate the convergence~\citep{franccois2015discount}. Studying such properties in classical LQR problems will be our important future work.

\bibliographystyle{plain}
\bibliography{mybibfile}

\appendix

\section{Policy gradient methods and sample complexity}


To make our presentation self-contained, we first provide a well-known fact~\citep{malik2019derivative,perdomo2021stabilizing} as a basis for the policy gradient analysis. With slight abuse of notation, we use $J_{\gamma}(K,A,B)$ to denote the discounted cost under system model $(A,B)$.
\begin{lemma}
Suppose that $\sqrt{\gamma}\rho(A-BK)<1$. Then, it holds $J_{\gamma}(K,A,B) = J_{1}(K, \sqrt{\gamma}A, \sqrt{\gamma}B)$.
\end{lemma}
\begin{proof}
	By assumption, $J_{\gamma}(K,A,B) < \infty$. Then, it follows from the definition that
	\begin{align*}
		J_{\gamma}(K,A,B) &= \mathbb{E}_{x_0} \sum_{t=0}^{\infty}\gamma^{t} (x_{t}^{\top} Q x_{t}+u_{t}^{\top} R u_{t}) \\
		& = \mathbb{E}_{x_0} \sum_{t=0}^{\infty}\gamma^{t} x_{t}^{\top} (Q + K^{\top}RK) x_{t} \\
		& = \mathbb{E}_{x_0} \sum_{t=0}^{\infty} ((\sqrt{\gamma}A+\sqrt{\gamma}BK)^tx_{t})^{\top} Q ((\sqrt{\gamma}A+\sqrt{\gamma}BK)^tx_{t}) \\
		& = J_{1}(K,A,B).
	\end{align*} 
\end{proof}

Hence, the policy gradient analysis for the discounted LQR can be preserved from that for standard LQR problems with damped systems. Before we establish the sample complexity in Lemma \ref{lem:gradient}, we prove Lemma \ref{lem:barJ} which suggests a lower bound of the threshold $\bar{J}$ for the policy descent in Algorithm \ref{alg:InfPG}.

\begin{proof}(of Lemma \ref{lem:barJ})
	Let $K^*_{\gamma_2}$ be the optimal policy that minimizes $J_{\gamma_2}(K)$. For $0<\gamma_1<\gamma_2\leq$, we have that $\sqrt{\gamma_1}\rho(A-BK^*_{\gamma_2})<\sqrt{\gamma_2}\rho(A-BK^*_{\gamma_2})<1$. Then, for $0<\gamma_1<\gamma_2\leq$ it holds that
	\begin{align*}
		J^*_{\gamma_1} & < J_{\gamma_1}(K^*_{\gamma_2}) \\
		 &\leq  \sum_{t=0}^{\infty}\gamma_1^{t} x_{t}^{\top} (Q + K^{*\top}_{\gamma_2}RK^*_{\gamma_2}) x_{t} \\
		& \leq  \sum_{t=0}^{\infty}\gamma_2^{t} x_{t}^{\top} (Q + K^{*\top}_{\gamma_2}RK^*_{\gamma_2}) x_{t}  \\
		& = J_{\gamma_2}^*
	\end{align*}
	The proof is completed.
\end{proof}

Since policy gradient methods can achieve any desired accuracy in the cost, it suffices to set $\bar{J} \geq J_1^*$ to ensure the threshold condition (\ref{equ:Jbar}) in Algorithm \ref{alg:InfPG}.

Our sample complexity result in Lemma \ref{lem:gradient} for gradient descent is a refinement of \citet[Theorem 1]{mohammadi2020linear} by neglecting the problem-dependent constants, which builds on the assumption that the distribution of the initial state has a bounded sub-Gaussian norm. We assume that $\mathcal{D}$ has bounded supports $\|x_0\|\leq d$, which is a special case of their assumption. Let $K_0$ be the initial policy and $a = J_{\gamma}(K_0)$. Consider the policy gradient method at $k$-th iteration
\begin{equation}\label{equ:gradient_des}
K_{k+1} = K_k - \eta \hat{\nabla} J_{\gamma}(K_k), ~~k=0,1,2,\dots,
\end{equation}
with $\hat{\nabla} J_{\gamma}(K_k)$ estimated by Algorithm \ref{alg:gradient}. The result \citep[Theorem 1]{mohammadi2020linear} can be restated as follows.
\begin{lemma}\label{lem:moh}
	Let the simulation time $\tau$ and number of trajectories $M$ in Algorithm \ref{alg:gradient} satisfy
	$\tau \geq \theta'(a)\log(1/\epsilon) $ and $ M \geq c(1+\beta^{4} d^{4} \theta(a) \log ^{6} n) n$ for a desired accuracy $\epsilon > 0$ and some $\beta >0$. Then,  for a smoothing radius $r<\theta''(a)\sqrt{\epsilon}$ and some constant step size $\eta = 1/(\omega(a)L(a))$, iterations (\ref{equ:pg}) achieves 
	$J_{\gamma}(K_{k}) - J_{\gamma}^* \leq \epsilon$
	in at most $$k \leq -\log(\epsilon^{-1}(J_{\gamma}(K_0)-J_{\gamma}^*))/\log(1-\mu(a)\eta/8)$$ iterations. This holds with probability not smaller than
	$1-c' k(n^{-\beta}+M^{-\beta}+M e^{-\frac{n}{8}}+e^{-c' M})$. Here, $\omega(a)=c''(\sqrt{m}+\beta d^{2} \theta(a) \sqrt{m n} \log n)^{2}$, $c,c',c''$ are constants, $\mu(a)$ and $L(a)$ are the gradient dominance and smoothness parameters of $J_{\gamma}(K)$ over the sublevel set $\{K | J_{\gamma}(K) \leq J_{\gamma}(K_0)\}$, and $\theta, \theta', \theta''$ are positive polynomials that depend only on the parameters of the discounted LQR problem.
\end{lemma}

Lemma \ref{lem:moh} significantly improves from the existing literature in both the simulation time and the number of sampled trajectories. Specifically, the total number of trajectories is only $\mathcal{O}(\log(1/\epsilon))$ for an accuracy $\epsilon$ compared to  $\mathcal{O}((1/\epsilon^4)\log(1/\epsilon))$ in \cite{fazel2018global} and $1/\epsilon$ in \cite{malik2019derivative}. Similarly, the required simulation time is only $\mathcal{O}(\log(1/\epsilon))$ in contrast to $\text{poly}(1/\epsilon)$ in \cite{fazel2018global} and infinite simulation time in \cite{malik2019derivative}. This breakthrough largely relies on the connection between the optimization landscape and a convex parameterization of the LQR problem. Base on it, \cite{mohammadi2021convergence} proposes to estimate only the direction of policy gradient instead of its value, which leads to much lower sample complexity.

\section{Extension to stochastic noise setting}
We now discuss the extension to the following linear dynamics with noises
\begin{equation}\label{equ:sys_noise}
	x_{t+1} = Ax_t + Bu_t + w_t, ~~x_0= 0,
\end{equation}
where $w_t \in \mathbb{R}^n$ is the additive noise. Here, we assume that the sequence $\{w_t\}$ is independently sampled from a distribution $\mathcal{D}_{\text{add}}$ with zero mean and unit variance. Also, its supports are bounded, i.e., $\|w_t\|\leq d$. In this case, the initial state is fixed to $x_0 = 0$.

Define the discounted LQR in the stochastic noise setting
$$
J_{\text{add}, \gamma}(K) := \mathbb{E}_{w} \sum_{t=0}^{\infty}\gamma^{t} (x_{t}^{\top} Q x_{t}+u_{t}^{\top} R u_{t})~~
\text {subject to} ~(\ref{equ:sys_noise}), x_0 = 0, w_t\sim \mathcal{D}_{\text{add}} , \text{and}~u_t = -K x_t.
$$
Due to the connections between the two settings \citep{malik2019derivative}, the closed-form of the discounted LQR cost can be provided below as a counterpart of Lemma \ref{lem:closed}.
\begin{lemma}
Suppose that $\sqrt{\gamma}\rho(A-BK)<1$. Then, the discounted LQR cost can be written as
\begin{equation}
J_{\text{add},\gamma}(K) = \frac{\gamma}{1-\gamma} \text{Tr}(P_K^{\gamma})
\end{equation}
where $P_K^{\gamma}$ is a unique positive definite solution to the Lyapunov equation 
$$
P_K^{\gamma} = Q + K^{\top}RK + \gamma (A-BK)^{\top}P_K^{\gamma}(A-BK).
$$
\end{lemma}

Clearly, $J_{\text{add},\gamma}(K)$ and $J_{\gamma}(K)$ are equal up to a coefficient. Hence, a counterpart of Theorem \ref{coro} can be derived for the stochastic noise setting.
\begin{theorem}
	Suppose that $\sqrt{\gamma}\rho(A-BK)<1$. Then, $J_{\text{add},\gamma'}(K)<+\infty$ if 
	\begin{equation}\label{equ:dis_gamma}
			\gamma' \leq (1 + \frac{\underline{\sigma}(Q+K^{\top}RK)}{(1/\gamma - 1)J_{\text{add},\gamma}(K) - \underline{\sigma}(Q+K^{\top}RK)}) \gamma.
	\end{equation}
\end{theorem}
\begin{proof}
	A sufficient condition for (\ref{equ:rule}) to hold is 
	$$
	1-\frac{\gamma}{\gamma'} < \frac{\underline{\sigma}(Q+K^{\top}RK)}{\|P\|},
	$$
	where $P = Q + K^{\top}RK + \gamma (A-BK)^{\top}P(A-BK)$. Since $\|P\|$ is upper bounded by
	$$
	\|P\| \leq \text{Tr}(P) = \frac{1-\gamma}{\gamma} J_{\text{add},\gamma}(K),
	$$
	The proof is completed.
\end{proof}
The discount policy gradient algorithm in this setting is provided in Algorithm \ref{alg:noise}. We note that the update rate (\ref{equ:rate_noise}) for the discount factor also has a lower bound. To see this, let $J_{\text{add},\gamma}(K) \leq \bar{J}$. Then, it follows that
$$
\alpha^i  \geq  \frac{\underline{\sigma}(Q)}{2(1/\gamma^0 - 1)\bar{J} - \underline{\sigma}(Q)},
$$
 which implies that the discount factor can be multiplied by a positive constant each iteration. Thus, a similar finite-time convergence guarantee as in Theorem \ref{the:finite-time} can be proved, which we omit here.

\begin{algorithm2e}[t]
	\caption{The discount policy gradient algorithm for the stochastic noise setting}
	\label{alg:noise}	
	\LinesNumbered
	\KwIn{Initial policy $K^0 = 0$ and discount factor $\gamma^0$, simulation time $\tau$, number of trajectories $N$ for cost evaluation.}
	\For {$i=0,1,\cdots$}
	{
		Evaluate $\hat{J}_{\text{add},\gamma^i}(K^{i}) = \frac{1}{N}\sum_{j=0}^{N-1} J_{\text{add},\gamma^i}^{\tau}(K^i, w^j) $ with a realization of noise sequence $w^j$ \;
		Compute the discount factor $\gamma^{i+1} = (1+\alpha^i)\gamma^i$ with $\alpha^i$ given by \
		\begin{equation}\label{equ:rate_noise}
			\alpha^i = \frac{\underline{\sigma}(Q+(K^{i})^{\top}RK^i)}{2(1/\gamma^i - 1)\hat{J}_{\text{add},\gamma^i}(K^{i}) - \underline{\sigma}(Q+(K^{i})^{\top}RK^i)};
		\end{equation}
		\If{$\gamma^{i+1} \geq  1 $}{Return a stabilizing policy $K^i$ \;}
		Update $K^{i+1}$ via policy gradient starting from policy $K^i$, such that 
		$$
		J_{\text{add},\gamma^{i+1}}( K^{i+1}) < \bar{J}.
		$$
	}
\end{algorithm2e}

We now analyze the function evaluation step in line 2 of Algorithm \ref{alg:noise} in the stochastic noise setting. For a realization of the noise sequence $w = \{w_k\}$, we define the truncated cost function
$$
J_{\text{add}, \gamma}^{\tau}(K,w) := \sum_{t=0}^{\tau -1}\gamma^{t} (x_{t}^{\top} Q x_{t}+u_{t}^{\top} R u_{t}),~~
\text {subject to} ~(\ref{equ:sys_noise}), x_0 = 0, \text{and}~u_t = -K x_t.
$$
We show that the bias $J_{\text{add}, \gamma}(K,w) - J_{\text{add}, \gamma}^{\tau}(K,w)$ also has exponential dependence in $\tau$ as in Lemma \ref{lem:error}. 
\begin{lemma}
	For a desired accuracy $0<\epsilon<1$, it holds that $J_{\text{add}, \gamma}(K,w) - J_{\text{add}, \gamma}^{\tau}(K,w) \leq \epsilon$ when the simulation time satisfies $\tau \geq \mathcal{O}(\log(1/\epsilon))$.
\end{lemma}
\begin{proof}
	Noting that $x_t = \sum_{k=1}^{t}(A-BK)^kw_{t-k}$, we have that for a policy with $\sqrt{\gamma}\rho(A-BK)<1$, 
	\begin{align*}
		&J_{\text{add},\gamma}(K, w) - J_{\text{add},\gamma}^{\tau}(K, w) \\
		&= \sum_{t=\tau}^{\infty}\gamma^{t} x_{t}^{\top} (Q+K^{\top}RK) x_{t} \\
		&\leq \|Q+K^{\top}RK\| \cdot \sum_{t=\tau}^{\infty}\gamma^{t}(\|\sum_{k=1}^{t}(A-BK)^kw_{t-k}\|^2) \\
		&\leq \|Q+K^{\top}RK\| d^2 \cdot \sum_{t=\tau}^{\infty} \gamma^{t}(\|\sum_{k=1}^{t}(A-BK)^k\|^2) \\
		&\leq \|Q+K^{\top}RK\| d^2 \cdot \sum_{t=\tau}^{\infty} \gamma^{t} \frac{\|A-BK\|^2(1-\|A-BK\|^{2t})}{1-\|A-BK\|^2} \\
		& = \|Q+K^{\top}RK\| d^2 \cdot \sum_{t=\tau}^{\infty}(\frac{ \gamma^{t} \|A-BK\|^2}{1-\|A-BK\|^2} - \frac{\|\sqrt{\gamma}(A-BK)\|^{2t}}{1-\|A-BK\|^2} ) \\
		& \leq  \frac{\|Q+K^{\top}RK\| d^2}{1-\|A-BK\|^2} \cdot (\frac{\gamma^{\tau} \|A-BK\|^2}{1-\gamma} - \frac{\|\sqrt{\gamma}(A-BK)\|^{2\tau}}{1-\|\sqrt{\gamma}(A-BK)\|^2}).
	\end{align*}

	By using the same techniques in the proof of Lemma \ref{lem:error}, $\|A-BK\|$ can be further upper bounded by a uniform constant, which we omit here for simplicity. Thus, the bias $J_{\text{add}, \gamma}(K,w) - J_{\text{add}, \gamma}^{\tau}(K,w)$ induced by finite simulation time decreases exponentially with respect to $\tau$. 
\end{proof}

One can further quantize the number of sampled trajectories $N$ in this setting to yield a similar result in Lemma \ref{lem:error} using concentration bounds. The proof details are omitted. 

Finally, we briefly discuss the analysis of policy gradient methods for Algorithm \ref{alg:noise}. Since \cite{mohammadi2021convergence, mohammadi2020linear} only consider the randomness of $x_0$, their convergence analysis might not be able to applied to the stochastic noise setting. Also, \cite{malik2019derivative} assumes that the simulation time is infinite to obtain a unbiased estimation for the cost, we cannot directly use their results. Nevertheless, one can adopt similar analysis techniques in \cite{fazel2018global} to provide a linear convergence rate for the gradient descent, though the sample complexity may be polynomial instead of linear in the desired accuracy.

\end{document}